\documentclass[12pt]{article}
\usepackage{amssymb}
\usepackage{mathrsfs}
\usepackage{pifont}
\usepackage{tipa}
\usepackage{amsmath}
\usepackage{amssymb,amsmath}
\usepackage{amsfonts}
\catcode`@=11
 \long\def\@makefntext#1{\noindent #1}

\begin{document}
\title{\large{\textbf{Unitary representations for the
 Schr\"{o}dinger-Virasoro\\ Lie algebra}} }
\author{Xiufu Zhang $^{1,2}$, Shaobin Tan $^1$\\
{\scriptsize 1. School of Mathematical Sciences, Xiamen University,
Xiamen 361005,  China}\\
{\scriptsize 2. School of Mathematical Sciences, Xuzhou Normal
University, Xuzhou 221116, China}}
\date{}
\maketitle \footnotetext{\footnotesize* Supported by the National
Natural Science Foundation of China (No. 10931006).}
\footnotetext{\footnotesize** Email: xfzhang@xznu.edu.cn;\quad
tans@xmu.edu.cn}

\numberwithin{equation}{section}

\begin{abstract}
In this paper, conjugate-linear anti-involutions and unitary
Harish-Chandra modules over the Schr\"{o}dinger-Virasoro algebra are
studied. It is proved that there are only two classes
conjugate-linear anti-involutions over the Schr\"{o}dinger-Virasoro
algebra. The main result of this paper is that a unitary
Harish-Chandra module over the Schr\"{o}dinger-Virasoro algebra is
simply a unitary Harish-Chandra module over the Virasoro algebra.
\vspace{2mm} \vspace{2mm}
\end{abstract}

\vskip 3mm \noindent{\section{{Introduction}}}

The Schr\"{o}dinger-Virasoro algebra $\mathfrak{sv}$ is defined to
be a Lie algebra with $\mathbb{C}$-basis
$\{L_n,M_n,Y_{n+\frac{1}{2}},c \mid n\in \mathbb{Z}\}$ subject to
the following Lie brackets:
\begin{align*}
&[L_m,L_n]=(n-m)L_{n+m}+\delta_{m+n,0}\frac{m^{3}-m}{12}c,\\
&[L_m,M_n]=nM_{n+m},\\
&[L_m,Y_{n+\frac{1}{2}}]=(n+\frac{1-m}{2})Y_{m+n+\frac{1}{2}},\\
&[Y_{m+\frac{1}{2}},Y_{n+\frac{1}{2}}]=(n-m)M_{m+n+1},\\
&[M_m,M_n]=0=[M_m,Y_{n+\frac{1}{2}}]=[\mathfrak{sv},c].
\end{align*}
It was introduced by M. Henkel in Ref. [7] by looking at the
invariance of the free Schr\"{o}dinger equation. Due to its
important roles in mathematics and statistical physics, it has
been studied extensively by many authors. In Refs. [17, 19, 20],
the twisted Schr\"{o}dinger-Virasoro algebra,
$\varepsilon$-deformation Schr\"{o}dinger-Virasoro algebra,  the
generalized Schr\"{o}dinger-Virasoro algebra and  the extended
Schr\"{o}dinger-Virasoro algebras are introduced. These Lie
algebras are all natural deformations of the
Schr\"{o}dinger-Virasoro algebra $\mathfrak{sv}.$ In Refs [6, 17,
19, 21], the derivations, the 2-cocycles, the central extensions
and the automorphisms for these algebras have been studied.

It is well known that the Virasoro Lie algebra is an important Lie
algebra, whose representation theory plays a crucial role in many
areas of mathematics and physics. Many representations, such as
the Harish-Chandra modules, the Verma modules and the Whittaker
modules, of it have been well studied (cf. Refs. [1, 3, 5, 8, 11,
13-16, 18]). The Virasoro Lie algebra is a subalgebra of the
Schr\"{o}dinger-Virasoro algebra, it is natural to consider those
representations for the Schr\"{o}dinger-Virasoro algebra. In Refs
[12, 19, 20, 22], the Harish-Chandra modules, the Verma modules,
the vertex algebra representations and the Whittaker modules over
the Schr\"{o}dinger-Virasoro algebra are studied.

In Ref. [1, 2-5], the nontrivial unitary irreducible unitary modules
are classified and the unitary highest weight modules of the
Virasoro algebra are well studied.  Motivated by these work, we
consider the unitary modules over the Schr\"{o}dinger-Virasoro
algebra in this paper. The paper is organized as follows. In section
2, we prove that there are only two classes conjugate-linear
anti-involutions over the Schr\"{o}dinger-Virasoro algebra
$\mathfrak{sv}$.  In section 3, we prove that a unitary weight
module over the Schr\"{o}dinger-Virasoro algebra is simply a unitary
weight module over the Virasoro algebra. Then the unitary weight
modules over $\mathfrak{sv}$ are classified since that the ones over
the Virasoro algebra are classified.

Throughout this paper we make a convention that the weight modules
over the Schr\"{o}dinger-Virasoro algebra and Virasoro algebra are
all with finite dimensional weight spaces, i.e., the Harish-Chandra
modules. The symbols $\mathbb{C}, \mathbb{N}$, $\mathbb{Z},$
$\mathbb{Z}_{+}$ and $\mathbb{Z}_{-}$ represent for the complex
field, the set of nonnegative integers, the set of integers, the set
of positive integers and the set of negative integers respectively.

\vskip 3mm\noindent{\section{Conjugate-linear anti-involution of
$\mathfrak{sv}$}}

It is easy to see the following facts about $\mathfrak{sv}:$

(i)  $\mathfrak{C}:=\mathbb{C}M_0\oplus\mathbb{C}c$ is the center of
$\mathfrak{sv}.$

(ii) If $x\in\mathfrak{sv}$ acts semisimply on $\mathfrak{sv}$ by
the adjoint action, then $x\in \mathfrak{h},$ where
$\mathfrak{h}:=\mathrm{span}_{\mathbb{C}}\{L_0,M_0,c\}$ is the
unique Cartan subalgebra of $\mathfrak{sv}$.

(iii) $\mathfrak{sv}$ has a weight space decomposition according to
the Cartan subalgebra $\mathfrak{h}:$
$$\mathfrak{sv}=\bigoplus_{n\in \mathbb{Z}}
\mathfrak{sv}_n\oplus\bigoplus_{n\in\mathbb{Z}}\mathfrak{sv}_{\frac{1}{2}+n},$$
where $\mathfrak{sv}_n=\mathrm{span}_{\mathbb{C}}\{L_n, M_n\},$
$\mathfrak{sv}_{\frac{1}{2}+n}=\mathrm{span}_{\mathbb{C}}\{Y_{\frac{1}{2}+n}\},$
$n\in \mathbb{Z}.$

If we denote $Vir=\oplus_{n\in \mathbb{Z}}\mathbb{C}L_n\oplus
\mathbb{C}c,$ $M=\oplus_{n\in \mathbb{Z}}\mathbb{C}M_n,
Y=\oplus_{n\in \mathbb{Z}}\mathbb{C}Y_{\frac{1}{2}+n}.$ Then we have
the following lemma:

\vskip 3mm\noindent{\bf{Lemma 2.1.}} $M\oplus Y\oplus \mathbb{C}c$
is the unique maximal ideal of $\mathfrak{sv}.$

\vskip 3mm\noindent{\bf{Proof.}} The proof is similar as that for
Lemma 2.2 in Ref. [19].\hfill$\Box$

\vskip 3mm\noindent{\bf{Definition 2.2.}} Let $\mathfrak{g}$ be a
Lie algebra and $\theta$ be a conjugate-linear anti-involution of
$\mathfrak{g}$, i.e. $\theta$ is a map
$\mathfrak{g}\rightarrow\mathfrak{g}$ such that
$$\theta(x+y)=\theta(x)+\theta(y),\ \theta(\alpha x)=\bar{\alpha}\theta(x),$$
$$\theta([x,y])=[\theta(y),\theta(x)],\ \theta^{2}=\mathrm{id}$$
for all $x,y\in \mathfrak{g},\ \alpha\in \mathbb{C},$ where
$\mathrm{id}$ is the identity map of $\mathfrak{g}.$ A module $V$
of $\mathfrak{g}$ is called unitary if there is a positive
definite Hermitian form $\langle\ ,\ \rangle$ on $V$ such that
$$\langle xu,v\rangle=\langle u,\theta(x)v\rangle$$ for all $u, v\in V, x\in \mathfrak{g}.$

\vskip 3mm\noindent{\bf{Lemma 2.3.}} (Proposition 3.2 in Ref. [1])
Any conjugate-linear anti-involution of $Vir$ is one of the
following types:

(i) $\theta_{\alpha}^{+}(L_n)=\alpha^{n}L_{-n},\
\theta_{\alpha}^{+}(c)=c,$ for some $\alpha\in \mathbb{R}^{\times},$
the set of nonzero real number.

(ii) $\theta_{\alpha}^{-}(L_n)=-\alpha^{n}L_{n},\
\theta_{\alpha}^{-}(c)=-c,$ for some $\alpha\in S^{1},$ the set of
complex number of modulus one.

\vskip 3mm\noindent{\bf{Lemma 2.4.}} Let $\theta$ be an arbitrary
conjugate-linear anti-involution of $\mathfrak{sv}.$ Then

(i) $\theta(M\oplus Y)=M\oplus Y.$

(ii) $\theta(\mathfrak{h})=\mathfrak{h}.$

(iii) $\theta(c)=\lambda c+\lambda^{'}M_0,$ $\theta(M_0)=\mu M_0,$
where $\lambda, \mu\in S^{1}, \lambda^{'}\in \mathbb{C}.$

\vskip 3mm\noindent{\bf{Proof.}} (i), $\forall x\in \mathfrak{sv},
y\in M\oplus Y,$ the identity $[x,\theta(y)]=\theta([y,\theta(x)])$
means that $\theta(M\oplus Y)$ is an ideal of $\mathfrak{sv}.$ Thus
$\theta(M\oplus Y)\subseteq M\oplus Y\oplus\mathbb{C}c$ by Lemma
2.1. Assume that
$$\theta(Y_{\frac{1}{2}+n})=a_n x+\beta_n c,$$ where $a_n,
\beta_n\in \mathbb{C}, x\in M\oplus Y.$ Then by $[L_0,
Y_{\frac{1}{2}+n}]=(\frac{1}{2}+n)Y_{\frac{1}{2}+n},$ we see that
$\beta_n=0.$ Moreover, $\theta(M_{n})\subseteq M\oplus Y$ since
$M=[Y,Y].$

For (ii), $\forall x\in \mathfrak{sv},$
$[x,\theta(M_0)]=\theta([M_0,\theta(x)])=0.$ So
$\theta(M_0)\in\mathfrak{C}.$ Similarly, $\theta(c)\in\mathfrak{C}.$
The identities $$[\theta(L_{0}),\theta(L_{n})]=-n\theta(L_{n}),\
[\theta(L_{0}),\theta(M_{n})]=-n\theta(M_{n}),$$ and
$$[\theta(L_{0}),\theta(Y_{\frac{1}{2}+n})]=-(\frac{1}{2}+n)\theta(Y_{\frac{1}{2}+n})$$
imply that $\theta(L_{0})$ acts semisimply on $\mathfrak{sv}.$ Thus
$\theta(L_{0})\in \mathfrak{h}.$

For (iii), note that $\mathfrak{C}$ is the center of
$\mathfrak{sv},$ we have $\theta(\mathfrak{C})=\mathfrak{C},$ so we
can assume
$$\theta(c)=\lambda c+\lambda^{'}M_0.$$ Since $\theta(M_0)\in
(M\oplus Y)\cap \mathfrak{C},$ we can assume $\theta(M_0)=\mu M_0.$
So $M_0=\theta^{2}(M_0)=\mu\bar{\mu}M_0,$ thus $\mu\in S^{1}.$
Similarly, we have $\lambda\in S^{1}.$\hfill$\Box$

\vskip 3mm\noindent{\bf{Proposition 2.5.}} Any conjugate-linear
anti-involution of $\mathfrak{sv}$ is one of the following types:
\begin{align*}(\mathrm{i}):\ \ &\theta_{\alpha,\beta,\mu}^{+}(L_n)
=\alpha^{n}L_{-n}+(\frac{n+1}{2}\alpha^{n-1}\beta
+\frac{n-1}{2}\alpha^{n-1}\mu\overline{\beta})M_{-n},\\
&\theta_{\alpha,\beta,\mu}^{+}(c)=c,\\
&\theta_{\alpha,\beta,\mu}^{+}(M_n)=\mu\alpha^{n}M_{-n},\\
&\theta_{\alpha,\beta,\mu}^{+}(Y_{\frac{1}{2}+n})=
\mu^{\frac{1}{2}}\alpha^{\frac{1}{2}+n}Y_{-\frac{1}{2}-n}
\end{align*}
for some $\alpha\in \mathbb{R}^{\times}, \mu\in S^{1},
\beta\in\mathbb{C}.$
\begin{align*}(\mathrm{ii}):\ \ &\theta_{\alpha,r_{1},r_{2},\mu}^{-}(L_n)
=-\alpha^{n}L_{n}+(\frac{n+1}{2}\alpha^{n+1}\mu r_{1}
-\frac{n-1}{2}\alpha^{n-1}\mu r_{2})M_{n},\\
&\theta_{\alpha,r_{1},r_{2},\mu}^{-}(c)=-c,\\
&\theta_{\alpha,r_{1},r_{2},\mu}^{-}(M_n)=\mu\alpha^{n}M_{n},\\
&\theta_{\alpha,r_{1},r_{2},\mu}^{-}(Y_{\frac{1}{2}+n})=
(-\mu)^{\frac{1}{2}}\alpha^{\frac{1}{2}+n}Y_{\frac{1}{2}+n}
\end{align*}
for some $\alpha, \mu\in S^{1}, r_{1},r_{2}\in\mathbb{R}.$

\vskip 3mm\noindent{\bf{Proof.}} Let $\theta$ be any
conjugate-linear anti-involution of $\mathfrak{sv}.$ By Lemma 2.4
(i), we have the induced conjugate-linear anti-involution of
$\mathfrak{sv}/(M\oplus Y)\simeq Vir:$
$$\bar{\theta}:\mathfrak{sv}/(M\oplus Y)\rightarrow\mathfrak{sv}/(M\oplus Y).$$
Thus by Lemma 2.3 we see that $\bar{\theta}$ is one of the following
types:

(a) $\bar{\theta}_{\alpha}^{+}(\bar{L_n})=\alpha^{n}\bar{L_{-n}},\
\bar{\theta}_{\alpha}^{+}(\bar{c})=\bar{c},$ for some $\alpha\in
R^{\times}.$

(b) $\bar{\theta}_{\alpha}^{-}(\bar{L_n})=-\alpha^{n}\bar{L_{n}},\
\bar{\theta}_{\alpha}^{-}(\bar{c})=-\bar{c},$ for some $\alpha\in
S^{1}.$

If $\bar{\theta}$ is of type (a),  we can assume
\begin{equation} \label{eq:1}
\theta(L_n)=\alpha^{n}L_{-n}+\sum_{i}\beta_{n,i}M_i+
\sum_{j}\gamma_{n,j} Y_{\frac{1}{2}+j},
\end{equation}
where $\beta_{n,i}, \gamma_{n,j}, a_n\in \mathbb{C}.$ By (2.1) and
Lemma 2.4 (ii), we have
$$\theta(L_0)=L_{0}+\beta_{0,0}M_0+a_0c.$$
Then by $[\theta(L_{-1}),\theta(L_{1})] =-2\theta(L_0),$ we deduce
that $a_0=0.$ Thus
\begin{equation} \label{eq:1}
\theta(L_0)=L_{0}+\beta_{0,0}M_0.
\end{equation}
By (2.1), (2.2) and the identity
$[\theta(L_n),\theta(L_0)]=n\theta(L_n),$ it can be deduced easily
that $\beta_{n,i}=0$ unless $i=-n,$ $\gamma_{n,j}=0$ for all $j\in
\mathbb{Z},$ i.e.,
\begin{equation} \label{eq:1}
\theta(L_n)=\alpha^{n}L_{-n}+\beta_{n,-n}M_{-n}.
\end{equation}
By (2.3) and the identity $[\theta(L_n),\theta(L_m)]
=(n-m)\theta(L_{m+n}) +\delta_{m+n,0}\frac{n-n^{3}}{12}\theta(c),$
we can get
\begin{eqnarray*}
((n-m)\beta_{m+n,-(m+n)}
-n\beta_{n,-n}\alpha^{m}+m\alpha^{n}\beta_{m,-m})M_{-(m+n)}
\end{eqnarray*}
\begin{equation} \label{eq:1}
=\delta_{m+n,0}\frac{n-n^{3}}{12}((1-\lambda)c-\lambda^{'}M_0)
\end{equation}
Let $m=-n\neq-1,0,1$ in (2.4), we see that
\begin{equation} \label{eq:1}
\lambda=1.
\end{equation}
Let $m=-n=1$ in (2.4), we have
\begin{equation} \label{eq:1}
\beta_{0,0}=\frac{\alpha^{-1}}{2}\beta_{1,-1}+\frac{\alpha}{2}\beta_{-1,1}.
\end{equation}
Let $m=-n=2$ in (2.4), we have
\begin{equation} \label{eq:1}
\lambda^{'}=8\beta_{0,0}-4\alpha^{2}\beta_{-2,2}-4\alpha^{-2}\beta_{2,-2}.
\end{equation}
Let $m=2, n=-1$ and $m=-2, n=1$ in (2.4) respectively , we have
\begin{equation} \label{eq:1}
\beta_{2,-2}=\frac{3\alpha}{2}\beta_{1,-1}-\frac{\alpha^3}{2}\beta_{-1,1},\
\
\beta_{-2,2}=\frac{3\alpha^{-1}}{2}\beta_{-1,1}-\frac{\alpha^{-3}}{2}\beta_{1,-1}.
\end{equation}
By (2.6)-(2.8), we have
\begin{equation} \label{eq:1}
\lambda^{'}=0.
\end{equation}
By (2.4), (2.5), (2.9) and Lemma 2.4 (iii), we have
\begin{equation} \label{eq:1}
\theta(c)=c,
\end{equation}
and
\begin{equation} \label{eq:1}
(n-m)\beta_{m+n,-(m+n)}
=n\beta_{n,-n}\alpha^{m}-m\alpha^{n}\beta_{m,-m}.
\end{equation}
Let $n=1$ in (2.11), we have
\begin{eqnarray*}
(1-m)\beta_{m+1,-(m+1)}+m\alpha\beta_{m,-m}-\alpha^{m}\beta_{1,-1}=0.
\end{eqnarray*}
Then using induction, we can prove that, for $m\geq1,$
\begin{equation} \label{eq:1}
\beta_{m,-m}=-(m-2)\alpha^{m-1}\beta_{1,-1}+(m-1)\alpha^{m-2}\beta_{2,-2}.
\end{equation}
Then  (2.8) and (2.12) give us  that
\begin{equation} \label{eq:1}
\beta_{m,-m}=\frac{m+1}{2}\alpha^{m-1}\beta_{1,-1}
-\frac{m-1}{2}\alpha^{m+1}\beta_{-1,1}, \ (\forall
m\in\mathbb{Z}_{+}).
\end{equation}
Let $n=-1$ in (2.11) and by a  similar argument as above, we can
prove that
\begin{equation} \label{eq:1}
\beta_{m,-m}=\frac{m+1}{2}\alpha^{m-1}\beta_{1,-1}
-\frac{m-1}{2}\alpha^{m+1}\beta_{-1,1}, (\forall
m\in\mathbb{Z}_{-}).
\end{equation}
Then by (2.6), (2.13) and (2.14), we see that
\begin{equation} \label{eq:1}
\beta_{m,-m}=\frac{m+1}{2}\alpha^{m-1}\beta_{1,-1}
-\frac{m-1}{2}\alpha^{m+1}\beta_{-1,1}, (\forall m\in\mathbb{Z}).
\end{equation}
Now by (2.3) and (2.15), we have
\begin{equation} \label{eq:1}
\theta(L_n)=\alpha^{n}L_{-n}+(\frac{n+1}{2}\alpha^{n-1}\beta_{1,-1}
-\frac{n-1}{2}\alpha^{n+1}\beta_{-1,1})M_{-n}.
\end{equation}

By Lemma 2.4 (i), we can assume that
\begin{equation} \label{eq:1}
\theta(M_m)=\sum_{i}\zeta_{m,i}M_i+\sum_{j}\xi_{m,j}Y_{\frac{1}{2}+j},
\end{equation}
\begin{equation} \label{eq:1}
\theta(Y_{\frac{1}{2}+n})=\sum_{i}\lambda_{n,i}M_i+
\sum_{j}\mu_{n,j}Y_{\frac{1}{2}+j},
\end{equation}
where $\zeta_{m,i},\ \xi_{m,j},\ \lambda_{n,i},\ \mu_{n,j}\in
\mathbb{C}.$ By (2.16), (2.17), Lemma 2.4 (iii) and the identity
$$[\theta(M_{-n}),\theta(L_{n})]
=-n\theta(M_{0})$$ we get that $\zeta_{-n,i}=0$ for all $i\neq n,$
$\alpha^{n}\zeta_{-n,n}=\mu,$ $\xi_{-n,j}=0$ for all $j\in
\mathbb{Z}.$ Thus
\begin{equation} \label{eq:1}
\theta(M_n)=\alpha^{n}\mu M_{-n}.
\end{equation}
By (2.3) and (2.19) we have $$L_1=\theta^{2}(L_1)
=L_1+(\alpha\beta_{-1,1}+\alpha^{-1}\mu\overline{\beta_{1,-1}})M_1.$$
Thus
\begin{equation} \label{eq:1}
\beta_{-1,1}=-\alpha^{-2}\mu\overline{\beta_{1,-1}}
\end{equation}
By (2.16) and (2.20), we have
\begin{equation} \label{eq:1}
\theta(L_n)=\alpha^{n}L_{-n}+(\frac{n+1}{2}\alpha^{n-1}\beta
+\frac{n-1}{2}\alpha^{n-1}\mu\overline{\beta})M_{-n},
\end{equation}
where $\beta=\beta_{1,-1}.$

By (2.16), (2.18) and the identity
$[\theta(Y_{\frac{1}{2}+m}),\theta(L_0)]
=(\frac{1}{2}+m)\theta(Y_{\frac{1}{2}+m}),$ we get that
$\lambda_{m,i}=0$ for all $i\in \mathbb{Z},$ $\mu_{m,j}=0$ unless
$j=-(m+1).$ Thus
\begin{equation} \label{eq:1}
\theta(Y_{\frac{1}{2}+m})=a_{\frac{1}{2}+m}Y_{-\frac{1}{2}-m},
\end{equation}
where $a_{\frac{1}{2}+m}=\mu_{m,-m-1}.$ If $n\neq1,$ by (2.16),
(2.22) and the identity $[\theta(Y_{\frac{1}{2}}),\theta(L_n)]
=(\frac{1-n}{2})\theta(Y_{\frac{1}{2}+n}),$ we have
$$\frac{1-n}{2}a_{\frac{1}{2}}\alpha^{n}Y_{-\frac{1}{2}-n}
=\frac{1-n}{2}a_{\frac{1}{2}+n}Y_{-\frac{1}{2}-n}.$$ Then
$a_{\frac{1}{2}+n}=\alpha^{n}a_{\frac{1}{2}},(n\neq1).$ By
$[\theta(Y_{\frac{3}{2}}),\theta(L_{-2})]
=\theta([L_{-2},Y_{\frac{3}{2}}]),$ we can easily get that
$a_{\frac{1}{2}+1}=\alpha a_{\frac{1}{2}}.$ So we have
\begin{equation} \label{eq:1}
a_{\frac{1}{2}+n}=\alpha^{n}a_{\frac{1}{2}}, \forall n\in
\mathbb{Z}.
\end{equation}
By (2.19), (2.22) and the identity
$[\theta(Y_{\frac{1}{2}+m}),\theta(Y_{-\frac{1}{2}-m})]
=\theta([Y_{-\frac{1}{2}-m},Y_{\frac{1}{2}+m}]),$ we have
$$a_{\frac{1}{2}+m}a_{-\frac{1}{2}-m}(2m+1)M_0=\mu(2m+1)M_0.$$ Thus
\begin{equation} \label{eq:1}
a_{\frac{1}{2}+m}a_{-\frac{1}{2}-m}=\mu.
\end{equation}
By (2.23) and (2.24), we see that
\begin{equation} \label{eq:1}
a_{\frac{1}{2}}=\sqrt{\mu\alpha}.
\end{equation}
By (2.22), (2.23) and (2.25), we have
\begin{equation} \label{eq:1}
\theta(Y_{\frac{1}{2}+m})
=\sqrt{\mu}\alpha^{\frac{1}{2}+m}Y_{-\frac{1}{2}-m}.
\end{equation}
Now (i) follows from (2.10), (2.19), (2.21) and (2.26).

\vskip 3mm If $\bar{\theta}$ is of type (b), by a similar discussion
in the way of (2.1)-(2.16), we can prove that
\begin{equation} \label{eq:1}
\theta(L_n) =-\alpha^{n}L_{n}+(\frac{n+1}{2}\alpha^{n-1}\beta_1
-\frac{n-1}{2}\alpha^{n+1}\beta_{-1})M_{n}.
\end{equation}
\begin{equation} \label{eq:1}
\theta(c)=-c,
\end{equation} where $\alpha\in S^{1}, \beta_1, \beta_{-1}\in \mathbb{C}.$
By a similar discussion in the way of (2.17)-(2.19) and
(2.22)-(2.26), we have
\begin{equation} \label{eq:1}
\theta(M_n)=\mu\alpha^{n}M_{n},
\end{equation}
\begin{equation} \label{eq:1}
\theta(Y_{\frac{1}{2}+n})=
(-\mu)^{\frac{1}{2}}\alpha^{\frac{1}{2}+n}Y_{\frac{1}{2}+n},
\end{equation}where $\mu\in S^{1}.$
By (2.27), (2.29) and the identities $\theta^{2}(L_{1})=L_1$ and
$\theta^{2}(L_{-1})=L_{-1}$, we see that
$$
\overline{\alpha}\beta_1=\overline{\beta_{1}}\alpha\mu,\
\alpha\beta_{-1}=\overline{\beta_{-1}}\overline{\alpha}\mu.
$$ If we set $\alpha=e^{i\sigma}, \mu=e^{i\tau}, \beta_1=|\beta_1|e^{ix},$
then by $\overline{\alpha}\beta_1=\overline{\beta_{1}}\alpha\mu,$ we
see that
$$\beta_1=|\beta_1|e^{i(2\sigma+\tau)}\ \mathrm{or}\ -|\beta_1|e^{i(2\sigma+\tau)}
.$$ Similarly, $$\beta_{-1}=|\beta_{-1}|e^{i(-2\sigma+\tau)}\
\mathrm{or}\ -|\beta_{-1}|e^{i(-2\sigma+\tau)}.$$ We set $r_1,
r_2\in \mathbb{R}$ such that $|r_1|=|\beta_{1}|$ and
$|r_2|=|\beta_{-1}|.$ Then
\begin{equation} \label{eq:1}
\beta_1=r_1\alpha^{2}\mu,\ \beta_{-1}=r_2\overline{\alpha}^{2}\mu.
\end{equation}
Thus (ii) follows from (2.27)-(2.31). \hfill$\Box$

\vskip 3mm The following Lemma is crucial for the proof of
Proposition 3.4.

\noindent{\bf{Lemma 2.6.}} Let $\theta$ be a conjugate-linear
anti-involution of the Schr\"{o}dinger-Virasoro algebra
$\mathfrak{sv}.$

(i) If $\theta=\theta_{\alpha,\beta,\mu}^{+},$  we denote by
$Vir^{'}$ the subalgebra of $\mathfrak{sv}$ generated by
$$\{c, L_{n}^{'}:=L_n-\frac{n-1}{2}\alpha^{-1}\beta M_{n}\mid
n\in\mathbb{Z}\}.$$ Then $Vir^{'}\simeq Vir$ and
$\theta_{\alpha,\beta,\mu}^{+}(L_{n}^{'})=\alpha^{n}L_{-n}^{'},
\theta_{\alpha,\beta,\mu}^{+}(c)=c.$

(ii) If $\theta=\theta_{\alpha,r_1,r_{2},\mu}^{-},$  we denote by
$Vir^{'}$ the subalgebra of $\mathfrak{sv}$ generated by
$$\{c, L_{n}^{'}:=L_n+x_nM_{n}\mid
n\in\mathbb{Z}\},$$ where $x_n\in \mathbb{C}$ satisfying
$\overline{x_{n}}\mu^{\frac{1}{2}}+x_n\mu^{-\frac{1}{2}}
=\frac{n-1}{2}r_2-\frac{n+1}{2}r_1.$ Then $Vir^{'}\simeq Vir$ and
$\theta_{\alpha,r_1,r_2,\mu}^{-}(L_{n}^{'})=-\alpha^{n}L_{n}^{'},
\theta_{\alpha,r_1,r_2,\mu}^{-}(c)=-c.$

\vskip 3mm\noindent{\bf{Proof.}} It can be checked directly, we omit
the details.\hfill$\Box$

\vskip 3mm\noindent{\bf{Lemma 2.7.}} (Proposition 3.4 in Ref. [1])
Let $V$ be a nontrivial irreducible weight $Vir$-module.

(i) If $V$ is unitary for some conjugate-linear anti-involution
$\theta$ of $Vir$,  then $\theta=\theta_{\alpha}^{+}$ for some
$\alpha>0.$

(ii) If $V$ is unitary for $\theta_{\alpha}^{+}$ for some
$\alpha>0,$ then $V$ is unitary for $\theta_{1}^{+}.$

\vskip 3mm\noindent{\bf{Proposition 2.8.}} Let $V$ be a nontrivial
irreducible weight $\mathfrak{sv}$-module.

(i) If $V$ is unitary for some conjugate-linear anti-involution
$\theta$ of $\mathfrak{sv}$, then
$\theta=\theta_{\alpha,\beta,\mu}^{+}$ for some $\alpha>0.$

(ii) If $V$ is unitary for $\theta_{\alpha,\beta,\mu}^{+}$ for some
$\alpha>0,$ then $V$ is unitary for $\theta_{1,\beta,\mu}^{+}.$

\noindent{\bf{Proof.}} (i) Suppose $V$ is unitary for some
conjugate-linear anti-involution $\theta$ of $\mathfrak{sv}$. By
Lemma 2.6, $V$ can be viewed as a unitary $Vir^{'}$-module for the
conjugate-linear anti-involution $\theta|_{Vir^{'}}$. Then $V$ is a
direct sum of irreducible unitary $Vir^{'}$-modules since any
unitary weight $Vir$-module is complete reducible. We claim that $V$
is a nontrivial $Vir^{'}$-module. Otherwise, for any $0\neq v\in V,
n\in\mathbb{Z}\setminus{0}, m\in\mathbb{Z},$ we have
\begin{align*}
&M_0v=-\frac{1}{n}(L_n^{'}M_{-n}-M_{-n}L_n^{'})v=0,\\
&M_nv=\frac{1}{n}(L_0^{'}M_n-M_nL_0^{'})v=0,\\
&Y_{\frac{1}{2}+m}v=\frac{2}{1+2m}(L_0^{'}Y_{\frac{1}{2}+m}
-Y_{\frac{1}{2}+m}L_0^{'})v=0.
\end{align*}So $\mathfrak{sv}.V=0,$  a contradiction. Thus there is a nontrivial
irreducible unitary $Vir^{'}$-submodule of $V$ for conjugate-linear
anti-involution $\theta|_{Vir^{'}}.$ By Lemma 2.7,
$\theta|_{Vir^{'}}=\theta_{\alpha}^{+}$ for some $\alpha>0.$ Then by
Proposition 2.5, we have $\theta=\theta_{\alpha,,\beta,\mu}^{+}$ for
some $\mu\in S^{1}, \beta\in\mathbb{C}.$

(ii) Suppose $V$ is unitary for $\theta_{\alpha,,\beta,\mu}^{+}$
for some $\alpha\in \mathbb{R}^{\times}, \mu\in S^{1},
\beta\in\mathbb{C}$ and $\langle\ .\ \rangle_{\alpha}$ is the
Hermitian form on $V.$ We can assume $V$ is generated by a
$L_0^{'}$-eigenvector $v_0$ with eigenvalue $a\in\mathbb{C}$ since
$V$ is irreducible weight $\mathfrak{sv}$-module. By
$$\langle L_0^{'}v_0,v_0\rangle_{\alpha}=\langle v_0,L_0^{'}v_0\rangle_{\alpha},$$ we see that $a\in
\mathbb{R}.$ Then $L_0^{'}$-eigenvalues on $V$ are of the form
$a+\frac{n}{2},$ $n\in\mathbb{Z}.$ Define a new form $\langle\ ,\
\rangle$ on $V$ by
$$\langle v,w\rangle=\alpha^{-\frac{n}{2}}\langle v,w\rangle_{\alpha},
\forall v, w\in \{v\mid L_0^{'}v=(a+\frac{n}{2})v\}.$$ It is easy
to check that this form makes $V$ unitary with the
conjugate-linear anti-involution $\theta_{1,\beta,\mu}^{+}.$
\hfill$\Box$

\vskip 5mm\noindent{\section{Unitary representations for
$\mathfrak{sv}$}}

In this section, we study the unitary weight modules for
$\mathfrak{sv}$. By Prop. 2.8, we see that the conjugate-linear
anti-involution is of the form $\theta_{\alpha,\beta,\mu}^{+}$ for
some $\alpha>0.$  For the sake of simplicity, we write
$\theta_{\alpha,\beta,\mu}^{+}$ by $\theta$.

\vskip 2mm It is known that a unitary weight module over Virasoro
algebra is completely reducible. This result also holds for
$\mathfrak{sv}:$

\vskip 2mm\noindent{\bf{Lemma 3.1.}} If $V$ is a unitary weight
module for $\mathfrak{sv}$, then $V$ is completely reducible.

\vskip 3mm\noindent{\bf{Proof.}} Let $N$ be a submodule. Then
$N^{+}:=\{v\in V|\langle v,N\rangle=0\}$ is a submodule of $V$
since for any $v\in N^{+},$ $\langle x.v,N\rangle=\langle
v,\theta(x)N\rangle=0.$ It is well known that any submodule of a
weight module is a weight module. For any weight $\lambda$ of $V,$
denote by $V_{\lambda}, N_{\lambda}, N^{+}_{\lambda}$ the weight
space with weight $\lambda$ of $V, N, N^{+}$ respectively.  It is
obvious that $dim(V_{\lambda})<\infty$ since $V$ is a
Harish-Chandra module, so we can extend an orthogonal basis of
$N_{\lambda}$ as an orthogonal basis of $V_{\lambda},$ thus we
have
$$V_{\lambda}=N_{\lambda}\oplus N^{+}_{\lambda},$$
which means $V=N\oplus N^{+}.$ \hfill$\Box$

\vskip 3mm\noindent{\bf{Lemma 3.2.}} (Theorem 1.3 (i) in [12]) An
irreducible weight module over $\mathfrak{sv}$ is either a
highest/lowest weight module or a uniformly bounded one.

It is well known (See Refs. [10] and [11]) that there are three
types modules of the intermediate series over $Vir$, denoted
respectively by $A_{a,b}, A_{\alpha}, B_{\beta},$ they all have
basis $\{v_k\mid k\in\mathbb{Z}\}$ such that $c$ acts trivially and
\begin{align*}
&A_{a,b}: L_nv_k=(a+k+nb)v_{n+k};\\
&A_{\alpha}: L_nv_k=(n+k)v_{n+k}\  \mathrm{if} \ k\neq0,\ L_nv_0=n(n+a)v_n;\\
&B_{\beta}: L_nv_k=kv_{n+k}\  \mathrm{if} \ k\neq-n,\
L_nv_{-n}=-n(n+a)v_0.
\end{align*}for all $n, k\in \mathbb{Z}.$
For the irreducible modules of the intermediate series of type
$A_{a,b},$ we have fact that: $A_{a,b}$ and $A_{c,d}$ are isomorphic
if and only if $a-c\in\mathbb{Z}$ and $b=d$ or $1-d.$

\vskip 3mm\noindent{\bf{Lemma 3.3.}} (Theorem 0.5 in [1]) Let $V$ be
an irreducible unitary module of $Vir$ with finite-dimensional
weight spaces. Then either $V$ is highest or Lowest weight, or $V$
is isomorphic to $A_{a,b}$ for some $a\in \mathbb{R}, b\in
\frac{1}{2}+\sqrt{-1}\mathbb{R}.$

\vskip 3mm\noindent{\bf{Proposition 3.4.}} A unitary weight module
over $\mathfrak{sv}$ is simply a unitary weight module over $Vir.$
That is, if $V$ is a unitary weight module over $\mathfrak{sv},$
then $M.V=Y.V=0.$

\vskip 2mm\noindent{\bf{Proof.}} Let $V$ be a unitary weight module
over $\mathfrak{sv}$ for a conjugate-linear anti-involution
$\theta.$ By Lemma 2.6 (i), $V$ is also unitary for $Vir^{'},$ thus
the well known result for the unitary modules over Virasoro Lie
algebra can be used freely. By Lemma 3.1, we may assume that $V$ is
irreducible. By Lemma 3.2, it is sufficient to consider the
following two cases:

\vskip 3mm\noindent{\bf{Case 1.}} $V$ is a unitary irreducible
highest/lowest weight module. Let $v_{\lambda}$ be a highest
weight vector. For $n\in \mathbb{Z}_{+},$ we have $\langle
M_{-n}v_{\lambda},M_{-n}v_{\lambda}\rangle =\langle
v_{\lambda},\mu M_{-n}M_{n}v_{\lambda}\rangle=0,$ thus
$M_{-n}v_{\lambda}=0.$ Furthermore,
$$\langle L_{-n}^{'}v_{\lambda},
M_{-n}v_{\lambda}\rangle=\langle
v_{\lambda},\alpha^{-n}L_{n}^{'}M_{-n}v_{\lambda}\rangle
=-n\alpha^{-n}\lambda(M_0)\langle
v_{\lambda},v_{\lambda}\rangle=0.$$ So $M_0v_{\lambda}=0.$ Thus
$$M.V=0.$$ For $n\in \mathbb{N},$ note that $M_0v_{\lambda}=0,$ we
have
$$\langle Y_{-\frac{1}{2}-n}v_{\lambda},Y_{-\frac{1}{2}-n}v_{\lambda}\rangle
=\langle
v_{\lambda},\mu^{\frac{1}{2}}Y_{\frac{1}{2}+n}Y_{-\frac{1}{2}-n}v_{\lambda}\rangle=0.$$
Thus $Y_{-\frac{1}{2}-n}v_{\lambda}=0,$ which means that
$$Y.V=0.$$

\vskip 3mm\noindent{\bf{Case 2.}} $V$ is a unitary irreducible
uniformly bounded module. As $Vir^{'}$-module, $V$ is a direct sum
of unitary irreducible $Vir^{'}$-submodules, So by Lemma 3.3 we
can suppose that
$$V=A_{a_1,b_1}\oplus\cdots \oplus A_{a_K,b_K}\oplus W,$$
where $a_i\in\mathbb{R}, b_i\in\frac{1}{2}+\sqrt{-1}\mathbb{R},$ $W$
is a trivial $Vir^{'}$-module.

\vskip 3mm\noindent{\bf{Subcase 2.1.}} $W=0.$

\vskip 3mm\noindent In this subcase
$$V=A_{a_1,b_1}\oplus\cdots \oplus A_{a_K,b_K}.$$
Let $\{v_k\mid k\in \mathbb{Z}\}$ be a basis of $A_{a_1,b_1}$ such
that $L_n^{'}v_k=(a_1+k+nb_1)v_{n+k}.$ As an irreducible
$\mathfrak{sv}$-module, $V$ is generated by the
$L_0^{'}$-eigenvector $v_0$ with eigenvalue $a_1.$ Thus
$L_0^{'}$-eigenvalue on $V$ are of the form $a_1+\frac{n}{2}, n\in
\mathbb{Z}.$ This means that $a_i\in \{a_1+\frac{n}{2}\mid n\in
\mathbb{Z}\},$ $i=1,\cdots,K.$ Recall that $A_{a+n,b}\simeq
A_{a,b}$ for any $n\in \mathbb{Z}.$  So there exists $0\leq
a<\frac{1}{2}$ such that $A_{a_i,b_i}$ are of the form $A_{a,b_i}$
or $A_{\frac{1}{2}+a,b_i}.$ i.e.,
$$V=A_{a,b_1}\oplus\cdots \oplus A_{a,b_R}\oplus
A_{\frac{1}{2}+a,d_1}\oplus\cdots \oplus A_{\frac{1}{2}+a,d_S},$$
where $0\leq a<\frac{1}{2}, b_i,
d_j\in\frac{1}{2}+\sqrt{-1}\mathbb{R}.$ Now we choose a basis
$$\{v_{k,l}\mid k\in\mathbb{Z}, 1\leq l\leq R\}\cup
\{v_{\frac{1}{2}+k,l^{'}}\mid k\in\mathbb{Z}, 1\leq l^{'}\leq S\}$$
such that
\begin{equation} \label{eq:1}
L_{m}^{'}v_{k,l}=(a+k+mb_l)v_{k+m,l},
\end{equation}
\begin{equation} \label{eq:1}
L_{m}^{'}v_{\frac{1}{2}+k,l^{'}}=(\frac{1}{2}+a+k+md_{l^{'}})v_{\frac{1}{2}+k+m,l^{'}}
\end{equation}
for $m\in\mathbb{Z}, 1\leq l\leq R, 1\leq l^{'}\leq S.$ Suppose
\begin{equation} \label{eq:1}
Y_{\frac{1}{2}}v_{k,l}=\sum_{l^{'}=1}^{S}\mu_{k,l}^{l^{'}}v_{\frac{1}{2}+k,l^{'}},
\end{equation}
\begin{equation} \label{eq:1}
Y_{\frac{1}{2}}v_{\frac{1}{2}+k,l^{'}}=\sum_{l=1}^{R}\lambda_{k,l^{'}}^{l}v_{k+1,l}.
\end{equation}
\vskip 3mm\noindent{\bf{Claim.}}
$Y_{\frac{1}{2}}v_{k,l}=0=Y_{\frac{1}{2}}v_{\frac{1}{2}+k,l^{'}},\
\forall k,l.$

Suppose the claim holds, Then $Y_{\frac{1}{2}}.V=0.$ Note that $Y,
M$ can be generated by $Y_{\frac{1}{2}}$ and $Vir^{'},$  we obtain
$M.V=0=Y.V,$ as desired. So it is sufficient to prove the claim.

\vskip 3mm\noindent{\bf{\emph{proof of the claim}.}} By (3.1), (3.3)
and the identity $[L_1^{'},Y_{\frac{1}{2}}]=0,$ we can easily deduce
that
\begin{equation} \label{eq:1}
(a+k+b_l)\mu_{k+1,l}^{l^{'}}=(a+\frac{1}{2}+k+d_{l^{'}})\mu_{k,l}^{l^{'}}.
\end{equation}
By (3.1)-(3,3) and the identity
$[L_{-1}^{'},Y_{\frac{1}{2}}]=Y_{-\frac{1}{2}},$ we have
\begin{equation} \label{eq:1}
Y_{-\frac{1}{2}}v_{k,l}=
\sum_{l^{'}=1}^{S}((a+\frac{1}{2}+k-d_{l^{'}})\mu_{k,l}^{l^{'}}-
(a+k-b_l)\mu_{k-1,l}^{l^{'}})v_{\frac{1}{2}+(k-1),l^{'}}.
\end{equation}
By (3.1), (3.2), (3.6) and the identity
$[L_1^{'},Y_{-\frac{1}{2}}]=-Y_{\frac{1}{2}},$ we have
\begin{eqnarray*}
&&\sum_{l^{'}=1}^{S}((a+\frac{1}{2}+k-d_{l^{'}})\mu_{k,l}^{l^{'}}
-(a+k-b_l)\mu_{k-1,l}^{l^{'}})
(a-\frac{1}{2}+k+d_{l^{'}})v_{\frac{1}{2}+k,l^{'}}-\\
&&(a+k+b_l)(\sum_{l^{'}=1}^{S}((a+\frac{1}{2}+k+1-d_{l^{'}})\mu_{k+1,l}^{l^{'}}-
(a+k+1-b_l)\mu_{k,l}^{l^{'}}))v_{\frac{1}{2}+k,l^{'}} \\
&&=-\sum_{l^{'}=1}^{S}\mu_{k,l}^{l^{'}}v_{\frac{1}{2}+k,l^{'}}.
\end{eqnarray*}
Thus
\begin{eqnarray*}
&&(a+\frac{1}{2}+k-d_{l^{'}})(a-\frac{1}{2}+k+d_{l^{'}})\mu_{k,l}^{l^{'}}
-(a+k-b_l)(a-\frac{1}{2}+k+d_{l^{'}})\mu_{k-1,l}^{l^{'}}\\
&&-(a+k+b_l)(a+\frac{1}{2}+k+1-d_{l^{'}})\mu_{k+1,l}^{l^{'}}
+(a+k+b_l)(a+k+1-b_l)\mu_{k,l}^{l^{'}}
\end{eqnarray*}
\begin{equation} \label{eq:1}
=-\mu_{k,l}^{l^{'}}
\end{equation}
By (3.1), (3.2), (3.4) and the identity
$[L_1^{'},Y_{\frac{1}{2}}]=0$ we get that
\begin{equation} \label{eq:1}
(a+\frac{1}{2}+k+d_{l^{'}})\lambda_{k+1,l^{'}}^{l}=(a+k+1+b_l)\lambda_{k,l^{'}}^{l},
\end{equation}
By (3.1)-(3.4) and the identity
$[L_{-1}^{'},Y_{\frac{1}{2}}]=Y_{-\frac{1}{2}}$ we have
\begin{equation} \label{eq:1}
Y_{-\frac{1}{2}}v_{\frac{1}{2}+k,l^{'}}
=\sum_{l=1}^{R}((a+k+1-b_{l})\lambda_{k,l^{'}}^{l}
-(\frac{1}{2}+a+k-d_{l^{'}})\lambda_{k-1,l^{'}}^{l})v_{k,l}.
\end{equation}
Then by  (3.1), (3.2), (3.4), (3.9) and identity
$[L_{1}^{'},Y_{-\frac{1}{2}}]=-Y_{\frac{1}{2}},$ we have
\begin{eqnarray*} &&(a+k+b_l)(a+k+1-b_l)\lambda_{k,l^{'}}^{l}
-(a+k+b_l)(a+\frac{1}{2}+k-d_{l^{'}})\lambda_{k-1,l^{'}}^{l}-\\
&&(a+k+\frac{1}{2}+d_{l^{'}})(a+k+2-b_{l})\lambda_{k+1,l^{'}}^{l}
+(a+k+\frac{1}{2}+d_{l^{'}})(a+\frac{3}{2}+k-d_{l^{'}})\lambda_{k,l^{'}}^{l}
\end{eqnarray*}
\begin{equation} \label{eq:1}
=-\lambda_{k,l^{'}}^{l}.
\end{equation}

If  $a-\frac{1}{2}+k+d_{l^{'}}\neq0$ for any $k\in\mathbb{Z},$  by
multiplying both sides of (3.7) by $a-\frac{1}{2}+k+d_{l^{'}}$ and
then using (3.5), we obtain that
$$
(-4k^{2}+\xi
k+\varsigma)\mu_{k,l}^{l^{'}}=-(a-\frac{1}{2}+k+d_{l^{'}})\mu_{k,l}^{l^{'}}
$$
for any $k\in\mathbb{Z},$ where $\xi, \varsigma\in \mathbb{C}.$ Thus
there exists at most two integers, say $k_1, k_2,$ such that for any
$k\in\mathbb{Z}\setminus\{k_1,k_2\},$ $\mu_{k,l}^{l^{'}}=0$ holds.
By (3.5), $\mu_{k_1,l}^{l^{'}}=0=\mu_{k_2,l}^{l^{'}}.$ Then
$$\mu_{k,l}^{l^{'}}=0, \forall k,l,l^{'}.$$ Thus by (3.3), we have
$$
Y_{\frac{1}{2}}v_{k,l}=0.
$$
By a similar discussion on $(3.2),(3.8)-(3.10),$ we get that
$$
Y_{\frac{1}{2}}v_{\frac{1}{2}+k,l^{'}}=0.
$$
\vskip 3mm\noindent If there exists $k\in\mathbb{Z}$ such that
$a-\frac{1}{2}+k+d_{l^{'}}=0,$  then we have
$$a=0, d_{l^{'}}=\frac{1}{2}$$ since $0\leq a<\frac{1}{2},
d_{l^{'}}\in\frac{1}{2}+\sqrt{-1}\mathbb{R}.$ Then by (3.5) and
(3.8) we have
\begin{equation} \label{eq:1}
\mu_{k,l}^{l^{'}}=0, \forall k\geq0.
\end{equation}
\begin{equation} \label{eq:1}
\lambda_{k,l^{'}}^{l}=0, \forall k\leq0.
\end{equation}
For $k=-1,$ note that $b_{l}\in\frac{1}{2}+\sqrt{-1}\mathbb{R},$ we
have
$$Y_{\frac{1}{2}}v_{-1,l}=-\frac{1}{b_l}Y_{\frac{1}{2}}L_{-1}^{'}v_{0,l}
=\frac{1}{b_l}Y_{-\frac{1}{2}}v_{0,l}
-\frac{1}{b_l}L_{-1}Y_{\frac{1}{2}}v_{0,l}
=\frac{1}{b_l}Y_{-\frac{1}{2}}v_{0,l},$$ Then by (3.3) and (3.6), we
have $$\mu_{-1,l}^{l^{'}}=0.$$ Note that $k-1+b_l,
-\frac{1}{2}+k+d_{l^{'}}\neq0$ for all $k\leq-1,$  so by (3.5)  we
have
\begin{equation} \label{eq:1}
\mu_{k,l}^{l^{'}}=0, \forall k<0.
\end{equation}
Since
$$Y_{\frac{1}{2}}v_{\frac{1}{2}+1,l^{'}}
=Y_{\frac{1}{2}}L_{1}^{'}v_{\frac{1}{2},l^{'}}
=L_{1}^{'}Y_{\frac{1}{2}}v_{\frac{1}{2},l^{'}}=0,$$ we have
$\lambda_{1,l^{'}}^{l}=0,$ then by (3.8) we have
\begin{equation} \label{eq:1}
\lambda_{k,l^{'}}^{l}=0, \forall k>0.
\end{equation}
From (3.11)-(3.14) we get that
$Y_{\frac{1}{2}}v_{k,l}=0=Y_{\frac{1}{2}}v_{\frac{1}{2}+k,l^{'}},\
\forall k,l,$ as required.

\vskip 3mm\noindent{\bf{Subcase 2.2.}} $W\neq0.$

\vskip 3mm\noindent Choose an arbitrary nonzero element $w\in W.$
$V$ generated by $w$ since  $V$ is an irreducible
$\mathfrak{sv}$-module. If $Y_{\frac{1}{2}}w=0,$ then $Y.W=M.W=0$
and $V$ is a trivial $\mathfrak{sv}$-module, a contradiction.  Thus
$Y_{\frac{1}{2}}w\neq0.$ By
$L_0^{'}Y_{\frac{1}{2}}w=\frac{1}{2}Y_{\frac{1}{2}}w,$ we see that
\begin{equation} \label{eq:1}
Y_{\frac{1}{2}}w\in A_{a_1,b_1}\oplus\cdots \oplus A_{a_K,b_K}.
\end{equation}
Moreover, $$A_{a_i,b_i}\simeq A_{0,b_i}\ \mathrm{or}\
A_{\frac{1}{2},b_i} $$ for each $i\in\{1,\cdots, K\}$ since $V$ is
generated by the eigenvector $w$ of $L^{'}_0$ with eigenvalue $0.$
So
$$V=A_{0,b_1}\oplus\cdots\oplus A_{0,b_R}\oplus A_{\frac{1}{2},d_1}\cdots\oplus A_{\frac{1}{2},d_S}\oplus W.$$
Choose the standard basis $\{v_{k,i}\mid k\in\mathbb{Z}\}$ and
$\{v_{\frac{1}{2}+k,j}\mid k\in\mathbb{Z}\}$ for each $A_{0,b_i}$
and $A_{\frac{1}{2},d_j}$ respectively. Suppose
$$
Y_{\frac{1}{2}}v_{k,l}=\sum_{l^{'}=1}^{S}\mu_{k,l}^{l^{'}}v_{\frac{1}{2}+k,l^{'}}+w_{k,l},
$$
$$
Y_{\frac{1}{2}}v_{\frac{1}{2}+k,l^{'}}=\sum_{l=1}^{R}\lambda_{k,l^{'}}^{l}v_{k+1,l}+w_{k,l^{'}}.
$$
where $w_{k,l}, w_{k,l^{'}}\in W.$ By a similar calculation as that
from identity (3.5) to identity (3.14) in Subcase 2.1 we have
\begin{equation} \label{eq:1}
w_{k,l}=0=Y_{\frac{1}{2}}v_{k,l},
\end{equation}
\begin{equation} \label{eq:1}
 w_{k,l^{'}}=0=Y_{\frac{1}{2}}v_{\frac{1}{2}+k,l^{'}}(k\neq
0),
\end{equation}
\begin{equation} \label{eq:1}
 w_{0,l^{'}}=Y_{\frac{1}{2}}v_{\frac{1}{2},l^{'}},
\end{equation}
and
\begin{equation} \label{eq:1}
Y_{-\frac{1}{2}}v_{\frac{1}{2}+k,l^{'}}=(k+1-d_{l^{'}})w_{k-1,l^{'}}.
\end{equation}
For any $m\in \mathbb{Z},$
$L_{m}^{'}Y_{\frac{1}{2}}v_{\frac{1}{2},l^{'}}=L_{m}^{'}w_{0,l^{'}}=0,$
so
$$Y_{\frac{1}{2}+m}v_{\frac{1}{2},l^{'}}=Y_{\frac{1}{2}}v_{\frac{1}{2}+m,l^{'}}=0$$
for $m\neq0,1.$ Thus
$$w_{0,l^{'}}=Y_{\frac{1}{2}}v_{\frac{1}{2},l^{'}}=-[L_{1},Y_{-\frac{1}{2}}]v_{\frac{1}{2},l^{'}}
=-L_{1}Y_{-\frac{1}{2}}v_{\frac{1}{2},l^{'}}+Y_{-\frac{1}{2}}L_{1}v_{\frac{1}{2},l^{'}}=(1+d_{l^{'}})(2-d_{l^{'}})w_{0,l^{'}}.$$
Note that $d_{l^{'}}\in\frac{1}{2}+\sqrt{-1}\mathbb{R},$ we have
\begin{equation} \label{eq:1}
w_{0,l^{'}}=0.
\end{equation}
By (3.16)-(3.20), we have
\begin{equation} \label{eq:1}
Y.(A_{0,b_1}\oplus\cdots\oplus A_{0,b_R}\oplus
A_{\frac{1}{2},d_1}\cdots\oplus A_{\frac{1}{2},d_S})=0,
\end{equation}
 and \begin{equation} \label{eq:1}
M.(A_{0,b_1}\oplus\cdots\oplus A_{0,b_R}\oplus
A_{\frac{1}{2},d_1}\cdots\oplus A_{\frac{1}{2},d_S})=0.
\end{equation}
By  (3.15),(3.21) and (3.22), we see that
\begin{equation} \label{eq:1}
Y(Y_{\frac{1}{2}}w)=M(Y_{\frac{1}{2}}w)=0.
\end{equation}
 Note that $Y_\frac{1}{2}w\neq0$ is also a
generator of $V$, combining with (3.23) and (3.15), we have
$$w\in U(L)(Y_{\frac{1}{2}}w)\subseteq A_{a_1,b_1}\oplus\cdots \oplus
A_{a_K,b_K},$$ this contradicts with that $0\neq w\in W$. Thus
Subcase 2.2 is impossible. This completes the proof of Proposition
3.4. \hfill$\Box$

\vskip 3mm If we denote the unitary weight modules over the
Schr\"{o}dinger-Virasoro algebra $\mathfrak{sv}$ by
$\overline{V}_{\lambda,0,0}, \underline{V}_{\lambda,0,0}$ and
$A_{a,b,0,0}$ corresponding respectively to the irreducible unitary
highest weight $Vir$-module $\overline{V}_{\lambda}$, the
irreducible unitary lowest weight $Vir$-module
$\underline{V}_{\lambda}$ and the irreducible unitary $Vir$-module
$A_{a,b}$. Then Proposition 3.4 and Lemma 3.3 give the
classification of the irreducible unitary weight modules over
$\mathfrak{sv}$:

\noindent{\bf{Theorem 3.5.}} An irreducible unitary weight module
$V$ over the Schr\"{o}dinger-Virasoro algebra is the highest weight
module $\overline{V}_{\lambda,0,0}$, or lowest weight module
$\underline{V}_{\lambda,0,0}$ for some $\lambda\in \mathfrak{h},$ or
$V$ is isomorphic to $A_{a,b,0,0}$ for some $a\in \mathbb{R}, b\in
\frac{1}{2}+\sqrt{-1}\mathbb{R}.$

\vskip 3mm \noindent{\section{{Note}}} This article has been
accepted by Journal of algebra and its applications.

\end{document}